\documentclass[12pt]{article}
\usepackage[parfill]{parskip}    

\usepackage{amssymb}

\newcommand{\ZZ}{{\mathbf{Z}}}
\newcommand{\CC}{{\mathbf{C}}}
\newcommand{\RR}{{\mathbf{R}}}
\newcommand{\proof}{\textbf{Proof\ }}
\newcommand{\proofoflem}{\textbf{Proof of lemma\ }}
\newcommand{\proofofcor}{\textbf{Proof of corollary\ }}
\newcommand{\qed}{\mbox{$\Box$}}

\newtheorem{thrm}{Theorem}[section]
\newtheorem{lem}[thrm]{Lemma}
\newtheorem{cor}[thrm]{Corollary}

\title{\bf Instanton homology and the Alexander polynomial}
\author{Yuhan Lim}
\date{Version: July 2009}                                        

\begin{document}

\maketitle

\section{Introduction}

In a recent paper \cite{KM}, Kronheimer and Mrowka revisit an instanton knot homology first defined and studied by A.\ Floer \cite{Fl}. For simplicity we assume the case of knots $K$ in the 3-sphere.
Kronheimer-Mrowka refine the instanton knot homology  theory by introducing a bigrading on the groups. The bigraded groups are denoted as $KHI(K)$. 
They conjecture that $KHI(K)$ is isomorphic to other knot homologies, such as the one defined by Ozsvath-Szarbo. A small step in this direction is to show that $KHI(K)$ recovers the Alexander polynomial.
In this note we prove this. In addition we state two simple corollaries of the result.

Concurrently, the main result of this note has also been obtained by Kronheimer and Mrowka in \cite{KM3}. 

\paragraph{Instanton homology for knots}
We review the version of instanton Floer homology  for an oriented knot $K$ in the 3-sphere
  $S^3$,  as considered in \cite{KM}.
It is defined to be the instanton homology of a certain  closure $K^T$ of the knot complement $S^3-\nu K^{\circ}$ by $({T-D^{\circ}})\times S^1$. Here $T$ is a 2-torus and $D$ is a small disk in $T$. 
Let $\Sigma_K$ be a Seifert surface for $K$, oriented so that the boundary, with the induced orientation, matches the orientation of $K$. In the closure, the boundaries of 
$\Sigma_K\cap(S^3-\nu K^{\circ})$ 
and $({T-D^{\circ}})\times\{1\}$ are identified as well as  meridians for $K$ on $\partial\nu(K)$ and the fibers $\{\mbox{\rm pt}\}\times S^1$ on the boundary of $({T-D^{\circ}})\times S^1$. An orientation of the 3-manifold is required in the definition of instanton homology. $K^T$ is given the orientation induced from $S^3-\nu K^{\circ}$, and $S^3$ is given the standard orientation on $\RR^3$, regarding it as $\RR^3\cup\{\infty\}$.

To define instanton homology in the non-homology 3-sphere case we also need to specify a complex line bundle $\omega\to K^T$ whose Chern class has a non-zero mod 2 evaluation on at least one integral homology class in $K^T$. Let $\alpha$ be a homologically non-trivial oriented simple closed curve in $T$. Choose $\omega_0$ to the complex line with Chern class Poincare dual to $\alpha$, thinking of $\alpha$ as living on $(T-D^{\circ})\times\{1\}$ in $K^T$. $I_*(K^T)_{\omega_0}$ denotes the instanton homology for $(K^T, \omega_0)$ with complex coefficients.

\paragraph{Generalized eigenspace decomposition}
There is only a relative mod 8 grading on $I_*(K^T)_{\omega_0}$ but there is an absolute mod 2 grading due to Froyshov \cite{F} (details below). We shall always assume the canonical mod 2 grading and any mod 8 grading used will be assumed consistent with the mod 2 grading. This makes the groups
$$
\widetilde{I}_0(K^T)_{\omega_0}=\bigoplus I_{2i}(K^T)_{\omega_0},
\quad
\widetilde{I}_1(K^T)_{\omega_0}=\bigoplus I_{2i-1}(K^T)_{\omega_0}
$$
well-defined.

Let $x_0$ be a point in $K^T$. The action of the $\mu$-map evaluated on $x_0$, $\mu(x_0)$
 sends $I_*(K^T)_{\omega_0}\to I_{*-4}(K^T)_{\omega_0}$.
By \cite{KM} this determines a splitting of  $\widetilde{I}_{*}(K^T)_{\omega_0}$
into $\pm 2$-eigenspaces, in their normalization.
Then by definition 
$$
KHI_*(K)=\widetilde{I}^+_{*}(K^T)_{\omega_0}=\mbox{+2-eigenspace.}
$$

Let $\widehat{\Sigma}$ denote the surface $(\Sigma_K\cap (S^3-\nu K^{\circ}))\cup (T-D^{\circ})\times\{1\}$, with the induced orientation from $\Sigma_K$.
The action of the $\mu$-map evaluated on $\widehat{\Sigma}$, $\mu(\widehat{\Sigma})$ sends
$I_*(K^T)_{\omega_0}\to I_{*-2}(K^T)_{\omega_0}$.
The action of $\mu(\widehat{\Sigma})$ and $\mu(x_0)$ commute; $\widetilde{I}^+_{*}(K^T)_{\omega_0}$ is preserved under the action of $\mu(\widehat{\Sigma})$. In \cite{KM} it is shown (with the normalization used there) that the eigenvalues of $\mu(\widehat{\Sigma})$ are the even integers $n$ satisfying the bound $|n|\leq 2g-2$, $g$ being the genus of $\widehat{\Sigma}$. Then there is a decomposition of $KHI_*(K)$ by the generalized eigenspaces of $\mu(\widehat{\Sigma})$, which in the notation of \cite{KM} is written
\begin{equation}\label{sum}
KHI _{*}(K)= \bigoplus_{i=-g}^{i=g} KHI_{*}(K,i).
\end{equation}  
Define the finite Laurent polynomial
$$
P_K (t) = P_{\omega_0}(K^T,\widehat{\Sigma})(t)= 
\chi_{-g} t^{-g} + \dots + \chi_0 + \chi_1 t +\dots + \chi_g t^g,
$$
where $\chi_i$ is the Euler characteristic of $ KHI_*(K,i)$.

\paragraph{Statement of results}

\begin{thrm}\label{mainthrm}
For any knot K in the 3-sphere, $P_K(t)$ is exactly the symmetrized and normalized Alexander polynomial $\Delta_K(t)$ of $K$.
\end{thrm}

This proves Conjecture 7.26 of \cite{KM}.
The symmetrized and normalized Alexander polynomial $\Delta_K(t)$ of $K$ is the unique representative that satisfies $\Delta_K(t^{-1})  =\Delta_K(t)$ and $\Delta_K(1)=1$. The proof of the theorem involves a straightforward application of Floer's surgery exact triangle.

Instead of the closure of the knot complement by $(T-D^{\circ})\times S^1$ we can also consider the standard closure by $D^2\times S^1$ given by 0-surgery on $K$, which we denote by $K^D$. We consider the instanton homology of $K^D$ but this time chose the complex line $\omega'$ with Chern class Poincare dual to an oriented meridian of $K$, regarding the meridian  as a curve in $K^D$.

\begin{cor}\label{maincor}
For oriented $K$, let $Q_K(t)$ be the finite Laurent polynomial defined analogously but with $K^D$ and complex line $\omega'$. Then 
$$
Q_K(t)=\frac{\Delta_K(t)-1}{t-2+t^{-1}}
$$
where $\Delta_K(t)$ is as above.
\end{cor}

For example the trefoil knot has $Q_K(t)=1$ and the figure-8 knot, $Q_K(t)=-1$.
We mention that for $I_*(K^D)_{\omega'}$ we do not actually know 
that $\mu({x_0})$ splits $I_*(K^D)_{\omega'}$ into $\pm 2$-eigenspaces. 
However it is not needed and we can simply use the generalized $+2$-eigenspace in place of the $+2$-eigenspace.

\begin{cor}
If the symmetrized and normalized Alexander polynomial $\Delta_K(t)$ of a knot K in $S^3$ is non-trivial i.e. $\Delta_K(t)\neq 1$, then $-1$-surgery on $K$ never yields a simply connected 3-manifold.
\end{cor}

\textbf{Remark } This corollary is a special case of the property P conjecture.
(Property P is proven in \cite{KM2} and \cite{KM}, independent of Perelman's proof of the Poincare conjecture. It uses a result of \cite{CG} that proves that consideration of $\pm 1$-surgery is sufficient.)

\proofofcor
If $\Delta_K(t)\neq 1$, then $Q_K(t)\neq 0$ and the instanton homology groups $I_*(K^D)_{\omega'}$ are non-trivial. By Floer's surgery exact triangle the instanton homology for the integral homology sphere $K_{-1}$ obtained by $-1$-surgery on $K$ is also non-trivial. Thus there must be at least one non-trivial representation of $\pi_1(K_{-1})$ into $SU(2)$. 
\qed

\section{Preliminaries}

\paragraph{Mod 2 grading} 

We briefly review the canonical mod 2 grading  since we will be using it throughout this note 
(\cite{F}, see also \cite[Sect.\ 6.5]{D}).
Let $[\varrho]\in I_*(Y)_{\omega}$. 
Suppose that $Y=\partial X$ as oriented manifolds and $E\to X$  a $U(2)$-bundle with connection $A$ that extends 
$\varrho$ on $Y$.

Let $\widehat{X}$ be the cylindrical-end manifold obtained by adjoining the semi-infinite tube $Y\times[0,\infty)$ to the boundary. Likewise extend $E$ to $\widehat{E}$ and also extend $A$ to $\widehat{A}$ by the pullback of $\varrho$ over the cylindrical end. We let $\mbox{Ind}{E}$ be the index of the anti-self dual operator on $\widehat{X}$ coupled to $\widehat{A}$. We also have an indicies 
$\mbox{\rm Ind}^{\pm}{X}$ of the anti-self dual operator on forms on $\widehat{X}$, on the positive/negative $\delta$-weighted spaces, where the weight is non-zero and sufficiently small in absolute value \cite[Sect.\ 3.3.1]{D}.

Define the mod 2 grading of $[\varrho]$ to be
\begin{displaymath}
\nu[\varrho] = \mbox{\rm Ind}{E} -3{\rm Ind}^{-}{X} \bmod 2.
\end{displaymath}
Index calculations (for instance \cite[Sect.\ 3.3.1]{D}) show that ${\rm Ind}^{-}{X}= b_1(X) -b_2^+(X)$. (Here we assume that $Y$ is connected.) $b_2^+$ is the dimension of a maximal positive definite subspace for the (possibly degenerate) intersection form on $H_2(X)$. 

Let us now suppose that $W$ is a cobordism between $Y$ and $Y'$ that induces a map
$I_W\colon I_{*}(Y)_{\omega}\to I_{*+k}(Y')_{\omega'}$. We wish to determine the value of $k \bmod 2$, the mod 2 degree of $I_W$.

\begin{lem} \label{degree}
The degree $k$ of the map $I_W$ above satisfies
$$
k = 3( b_1(W)-b_1(Y) +b_0(Y')-b_0(W)-b^+_2(W) )\bmod 2.
$$
\end{lem}

\textbf{Remark} If $Y$ or $Y'$ is disconnected the lemma is still valid. We need to however interpret  the instanton homology of a disjoint union $I_{*}(Y_0\cup Y_1)_{\omega\omega'}$ as the tensor product $I_{*}(Y_0)_{\omega}\otimes I_*(Y_1)_{\omega'}$. In particular the grading satisfies $\nu([\varrho]\otimes[\varrho'])=\nu[\varrho]+\nu[\varrho'] \bmod 2$.

\proofoflem
Let $[\varrho]\in I_{*}(Y)_{\omega}$ and $I_W([\varrho])=[\varrho']\in I_{*+k}(Y')_{\omega'}$.
Let $E_W$ be the $U(2)$-bundle with connection that limits to $\varrho$ and $\varrho'$ at the ends. Then  by assumption, ${\rm Ind}E_W=0$. Let $Y=\partial X$  and $E$ as above.  The additivity property of the index \cite[Sect.\ 3.3]{D} tells us that
\begin{eqnarray*}
{\rm Ind}(E\cup E_W) &=& {\rm Ind}E+{\rm Ind}E_W\\ 
{\rm Ind}^-(X\cup W) &=& {\rm Ind}^+X+{\rm Ind}^{-}W.
\end{eqnarray*}
On the other hand according to \cite[Prop.\ 3.10]{D},
\begin{displaymath}
{\rm Ind}^+X - {\rm Ind}^-X=-(b_0(\partial X)+b_1(\partial X))
\end{displaymath}
and  \cite[Prop.\ 3.15]{D} (with additional terms added for non-connected boundary)
\begin{displaymath}
{\rm Ind}^{-}W= b_1(W)-(b_0(W)-b_0(\partial W))-b_2^+(W).
\end{displaymath}
Since the difference $\nu[\varrho]-\nu[\varrho'] $ is given by
\begin{displaymath}
{\rm Ind}(E\cup E_W) -{\rm Ind}E -3({\rm Ind}^-(X\cup W) -{\rm Ind}^-X)\bmod 2,
\end{displaymath}
the result follows. \qed

\paragraph{Instanton invariants for 2-component links}

We will  need to introduce  versions of instanton homology associated to an oriented 2-component link $L$ in the 3-sphere. The definitions are parallel to the one for knots. These are defined as the instanton homology of certain closures of the link complement. Let $C=[0,1]\times S^1$ and $\Sigma_2$ the surface of genus 2. $D_1$ and $D_2$ are two small disjoint disks in $\Sigma_2$.
Let $\Sigma_L$ be an oriented (connected) Seifert surface with induced orientation on the boundary equal to the orientation on $L$. 
Set
\begin{eqnarray*}
L^C &=& (S^3-\nu L^{\circ})\cup (C\times S^1)\\
L^{\Sigma_2} &=& (S^3-\nu L^{\circ})\cup (\Sigma -D_1^{\circ} -D_2^{\circ})\times S^1.
\end{eqnarray*}
Again the union is taken along the boundaries such that  oriented boundary of
$\Sigma_L\cap (S^3-\nu L^{\circ})$
matches up with the oriented boundary of $[0,1]\times\{1\}$ or $(\Sigma -D_1^{\circ} -D_2^{\circ})\times\{1\}$. Additionally the meridians of the link should match up with the fibers $\{\mbox{\rm pt}\}\times S^1$ on the boundaries. In the case of $L^C$  the genus of $\widehat{\Sigma}$, the extension of $\Sigma_L$, is one greater that of $\Sigma_K$ but in the case of $L^{\Sigma_2}$ it is two greater.

As in the case for knots we can consider the instanton homology groups $I^+_{*}(L^C)_{\omega}$ and $I^+_*(L^{\Sigma_2})_{\omega}$ as the basis for the definition of instanton link homology. We shall define our choices $\omega=\omega_1, \omega_2$ below. Using the action of $\mu(\widehat{\Sigma})$ to give a decomposition  
$\widetilde{I}^+_{*}(L^C)_{\omega_1}$ and $\widetilde{I}^+_{*}(L^{\Sigma_2})_{\omega_2}$, we again can define the finite Laurent polynomials
\begin{displaymath}
P_{\omega_1}(L^C,\widehat{\Sigma})(t)\quad{\rm and}\quad  P_{\omega_2}(L^{\Sigma_2},\widehat{\Sigma})(t).
\end{displaymath}
Let $\alpha_0$ be an oriented simple arc in $S^3-\nu L^{\circ}$  connecting the two boundary components.
Let $\alpha_1$ be an oriented simple arc in $C\times S^1$ connecting the two boundary components of such that the boundary of $\alpha_1$ matches up with the boundary of $\alpha_0$ and the union is an oriented simple closed curve $\alpha'$. Then in the case of $L^C$ choose $\omega=\omega_1$ where $\omega_1$ has Chern class Poincare dual to $\alpha'$.

Think of $\Sigma_2-D_1^{\circ}-D_2^{\circ}$ as the connected sum of surfaces $C\sharp T$. Let $\alpha''$ be an oriented simple closed curve of the form $\alpha_0\cup \alpha_1$ where $\alpha_1$ is the arc as before but lives in the $C$ factor of $C\sharp T$. 
Let $\alpha$ be a homologically non-trivial simple closed curve in the $T$ factor of the connected sum. We think of $\alpha$ as living on $(\Sigma_2-D_1^{\circ}-D_2^{\circ})\times\{1\}$.
Then for $L^{\Sigma_2}$ choose $\omega=\omega_2$ where $\omega_2$ has Chern class Poincare dual to $\alpha'' + \alpha$.

\paragraph{The excision principle and instanton homology for $L^{\Sigma_2}$}
With our choices for $\omega$ it is actually the case that $P_{\omega_2}(L^{\Sigma_2},\widehat{\Sigma})(t)$ is determined by $P_{\omega_1}(L^C,\widehat{\Sigma})(t)$. This is a consequence of Floer's excision principle. 

\begin{lem}\label{product}
The finite Laurent polynomial 
$$
P_{\omega_2}(L^{\Sigma_2},\widehat{\Sigma})(t)=P_{\omega_1}(L^C,\widehat{\Sigma})(t)\:
\cdot P_{\omega_3}(\Sigma_2\times S^1,\Sigma)(t),
$$ 
where $P_{\omega_3}(\Sigma_2\times S^1,\Sigma_2)(t)$ is the finite Laurent polynomial for $\Sigma_2\times S^1$ derived from the instanton homology of $\Sigma_2\times S^1$ and defined in an analogous manner. The complex line 
$\omega_3$ has  Chern class the Poincare dual of the curve $\alpha$ above.
\end{lem}

\proof  Let $F_1, F_1$ be the two tori that form the boundary of $(S^3-\nu L^{\circ})\subset L^{\Sigma_2}$. Since $\omega_2$ evaluates non-trivally mod 2 on $F_i$ we may apply Floer's excision principle (see \cite[Theorem 7.7]{KM}) to conclude that there is an isomorphism
\begin{displaymath}
\widetilde{I}^+_{*}(L^{\Sigma_2})_{\omega_2}\cong
\widetilde{I}^+_{*}(\Sigma_2\times S^1)_{\omega_3}\otimes 
\widetilde{I}^{+}_{*}(L^C)_{\omega_1},
\end{displaymath}
However it is not immediately evident that this isomorphism preserves the mod 2 grading, we shall prove it below but only in this specific situation.

In the above isomorphism the action of $\mu(\widehat{\Sigma})$ on 
$\widetilde{I}^+_{*}(K^{\Sigma_2})_{\omega_2}$ corresponds to the action of 
$\mu(\widehat{\Sigma})\otimes 1+ 1\otimes \mu(\Sigma_2)$ on the tensor product space. It follows that the generalized $\lambda$-eigenspaces $W_{\lambda}$ in $\widetilde{I}^+_{*}(L^{\Sigma_2})_{\omega_2}$
 obey a relation of the form
\begin{displaymath}
W_{\lambda} \cong 
\bigoplus_{\lambda=\lambda_0+\lambda_1} U_{\lambda_0}\otimes V_{\lambda_1},
\end{displaymath}
where $U_{\lambda_0}$ and $V_{\lambda_1}$ are corresponding generalized eigenspaces in 
$\widetilde{I}^+_{*}(\Sigma_2\times S^1)_{\omega_3}$ and $\widetilde{I}^{+}_{*}(L^C)_{\omega_1}$ respectively.
The lemma follows easily, assuming the mod 2 grading claim. 

To establish the mod 2 grading, consider the surgery cobordism $W$ between $L^{C}\cup (\Sigma_2\times S^1)$ and $L^{\Sigma_2}$ in the proof of the excision principle. 
Let $H=T\times [-1,1]\times[-1,1]$. The boundary is 
$T\times\{-1,1\}\times[-1,1]\cup T\times[-1,1]\times\{-1,1\}$. We regard $T\times\{0\}\times[-1,1]$ as the \lq core\rq\ of $H$. Then $W$ is obtained from $(L^{C}\cup (\Sigma_2\times S^1))\times[0,1]$ by identifying $\nu F_1\cup \nu F_2$ in $(L^{C}\cup (\Sigma_2\times S^1))\times\{1\}$
with $T\times[-1,1]\times\{-1,1\}$ in $H$. Clearly $W$ deformation retracts onto $W_0$, the union of $L^{C}\cup (\Sigma_2\times S^1)$ and the core of $H$, where $F_1, F_2$ are identified with $\partial (T\times\{0\}\times[-1,1])$. 

In $L^C$ let $a$ be an oriented  simple closed curve corresponding to $S^1\times\{0\}\times\{1\}\subset C\times S^1$. Let $b$ be  an oriented simple closed curve corresponding to $\{\mbox{\rm pt}\}\times S^1\subset C\times S^1$. Then thinking of $a$ and $b$ as homology classes, we have $a=0$ in $H_1(L^C)$ and $b$ a generator of $H_1(L^C)$. On the other hand by the identification via the core of $H$, $a$ is a generator of $H_1(\Sigma_2\times\{\rm pt\})$ and $b$ corresponds to a fiber  $\{\mbox{\rm pt}\}\times S^1$. 

Since $H_1(L^C)\cong \ZZ\oplus\ZZ$ we have the rank of 
$H_1(\Sigma_2\times S^1)\oplus H_1(L^C)$ to be 7. Then  $H_1(W)=H_1(W_0)$ is obtained from 
$H_1(\Sigma_2\times S^1)\oplus H_1(L^C)$ by introducing the two relations above. Thus
$H_1(W_0)$ has rank 5. It is easily seen that $b^+_2(W)=0$, thus by Lemma~\ref{degree} the dimension shift is
\begin{eqnarray*}
k&=&b_1(W)-b_1(Y) +b_0(Y')-b_0(W)-b^+_2(W)\\
&=&5-(2+5)+1-1-0 \equiv 0\bmod 2.
\end{eqnarray*}
This completes the proof.
 \qed

\paragraph{Some calculations}

Lemma~\ref{product} above necessitates evaluation of $P_{\omega_3}(\Sigma_2\times S^1,\Sigma_2)(t)$. For late use we will also need to evaluate $P_{\omega_4}(T\times S^1)(t)$ where $\omega_4$ has Chern class dual to any homologically non-trivial oriented simple closed curve in $T$.

\begin{lem}\label{toruscal}
$P_{\omega_4}(T\times S^1)(t)=1$ for either orientation on $T\times S^1$.
\end{lem}

\proof After a diffeomorphism, the class $\omega_4$ can be assumed to have Chern class Poincare dual to the fibre $\{\mbox{\rm pt}\}\times S^1$.  It is well-known (for instance \cite[Prop.\ 1.14]{BD}) that the $SO(3)$-bundle with 2nd Steifel-Whitney class
$w_2\equiv c_1(\omega_4)\bmod 2$ carries a unique flat connection, up to gauge equivalence. In terms of $U(2)$-bundles this gives two flat connections $\varrho_0$ and $\varrho_1$ on the adjoint bundle. The instanton homology $I_{*}(T\times S^1)_{\omega_4}\cong \CC\oplus \CC$, where the dimensions differ by 4. The rest of the lemma follows easily once we can show that the generators lie in even dimensions.

Let $X=T\times D^2$ so that $Y=\partial X$. Let $E_{\varrho_i}\to X$ be a $U(2)$-bundle carrying connection $A$ that extends $\varrho_i$. Then ${\rm Ind}E_{\varrho_1}={\rm Ind}E_{\varrho_0}+4\bmod 8$ since the Floer dimensions differ by 4. We wish to evaluate ${\rm Ind}E_{\varrho_0}$. 
Let $\varphi\colon T\times S^1\to T\times S^1$ be the orientation reversing diffeomorphism that reverses the $S^1$ factor; denote by $\tilde{\varphi}$ a lift to the bundle level.
Then $\tilde{\varphi}^*(\varrho_0)$ is equivalent in instanton homology to either $\varrho_0$ or $\varrho_1$. Now form the double of $X$, identifying the boundaries via $\varphi$ and let $E=E_{\varrho_0}\cup E_{\tilde{\varphi}^*(\varrho_0)}$ be the corresponding bundle over the double. There are two ways of forming $E$, we choose the one that has $c_1$ dual to the surface $\{y_0\}\times D^2\cup \{y_0\}\times D^2$ in the double. Then
\begin{displaymath}
\mbox{\rm Ind}E=2(4c_2-c_1^2)-3(1-b_1+b^+_2) \equiv  0\bmod 8.
\end{displaymath}
 On the other hand, by the additivity of the index
 \begin{displaymath}
\mbox{\rm Ind}E_{\varrho_0} +\mbox{\rm Ind}E_{\tilde{\varphi}^*(\varrho_0)} = 
\mbox{\rm Ind}E_{\varrho_0} +\mbox{\rm Ind}E_{\varrho_0}+4n
\equiv 0\bmod 8.
\end{displaymath}
It follows that $\mbox{\rm Ind}E_{\varrho_0}\equiv 0\bmod 2$. Clearly $b_1(X)-b^+_2(X)=2$ so $\nu[\varrho_0]\equiv 0\bmod 2$.  \qed

\begin{lem}\label{sigma2}
$P_{\omega_3}(\Sigma_2\times S^1,\Sigma)(t) = t^{-1} -2 +t$ for either orientation on $\Sigma_2\times S^1$.
\end{lem}

\proof
According to \cite[Prop.\ 1.15]{BD} we have, in increasing order of dimensions
\begin{displaymath}
I_*(\Sigma_2\times S^1)_{\omega_3}
=\CC\oplus\CC^2\oplus\CC\oplus\{0\}\oplus\CC\oplus\CC^2\oplus\CC\oplus\{0\}.
\end{displaymath}
We shall fix this grading, beginning with zero for the first group reading left to right.  We assume for the moment that this is consistent with the mod 2 grading. We shall prove it below by showing that the group $\CC^2$ above must be an odd graded group.

The action of $\mu(\Sigma_2)$ has eigenvalues either $-2$, $0$ or $+2$ on 
$\widetilde{I}^+_*(\Sigma_2\times S^1)_{\omega_3}$. On the odd dimensions $\mu(\Sigma_2)$ necessarily shifts each summand of the odd groups
$\CC^2\oplus\{0\} \oplus\CC^2\oplus\{0\}$ to the next. Thus $\mu(\Sigma_2)$ is zero on
$\widetilde{I}^+_{1}(\Sigma_2\times S^1)_{\omega_3}$ and this comprises it's entire generalized eigenspace decomposition. On the other hand the $\pm 2$-eigenvalues for $\mu(\Sigma_2)$ are simple on $\widetilde{I}^+_{*}(\Sigma_2\times S^1)$ \cite[Prop.~7.4]{KM}, so only the $\pm 2$-eigenvalues on $\widetilde{I}^+_0(\Sigma_2\times S^1)_{\omega_3}$ have non-trivial eigenspaces, each of dimension 1.
This completely decomposes $\widetilde{I}^+_{0}(\Sigma_2\times S^1)_{\omega_3}$ and the result follows, modulo the claim regarding the mod 2 grading.

Let $\varrho$ be the (perturbed) flat connection that is a generator $I_1(\Sigma_2\times S^1)_{\omega_3}$. Let $\varphi\colon\Sigma\times S^1\to \Sigma_2\times S^1$ be the orientation reversing diffeomorphism that reverses the $S^1$ factor.
 Let $X$ be such that $\partial X=Y$ and as $E_{\varrho}\to X$ the bundle with connection that extends $\varrho$. Let  and $\tilde{\varphi}$ be a lift to a bundle map. It must be the case that 
$\tilde{\varphi}^*(\varrho)$ is a generator in dimension 1 or 5, so by the same reasoning in the preceding lemma, ${\rm Ind}E_{\tilde{\varphi}^*(\varrho)}\equiv {\rm Ind}E_{\varrho}\bmod 4$. Repeating the  doubling argument we find
\begin{equation}\label{indexE}
2\mbox{Ind}E_{\varrho} \equiv 2(4c_2-c_1^2)-3(1-b_1+b^+_2)\bmod 4,
\end{equation}
where the right-hand terms refer to the doubled manifold (and bundle). We now need to construct an appropriate $X$. This is done as follows.

Consider two copies of $T\times S^1$; denote these as $Y_0$ and $Y_1$. Let 
$S_i=\{z_0\}\times S^1$ be the fibers over $z_0$ in $T$ in $Y_i$. Consider $H=S^1\times D^2\times [-1,1]$. Clearly $\partial H=S^1\times S^1\times[-1,1]\cup S^1\times D^2\times\{-1,1\}$. Attach $H$ to the product cobordism $(Y_0\cup Y_1)\times [0,1]$ along $(Y_0\cup Y_1) \times\{1\}$ by identifying  
$\nu(S_0)\cup \nu(S_1)$ with $S^1\times D^2\times \{0,1\}$. This should be done so that the boundary of the \lq core\rq\ of $H$, $S^1\times\{0\}\times[-1,1]$ goes to $S_0\cup S_1$ in such a way that the orientations of $S_i$ match up. It will then be seen that the resultant 4-manifold is a cobordism $W$ between $Y_0\cup Y_1$ and $\Sigma_2\times S^1$. The $U(2)$-bundle over $\Sigma_2\times S^1$ with determinant line $\omega_3$ extends over $W$; this complex line must be trivial over say $Y_1$ and consequently equal to $\omega_4$ over $Y_0$.

We may fill in the boundary component $Y_1$  by $X_1=T\times D^2$ in the standard manner; however for $Y_0$ by the same space $X_0$ only after an appropriate diffeomorphism on the boundary (see proof of Lemma~\ref{toruscal}). This is required in order that the $U(2)$-bundle over $W$ will in turn extend over the $X_i$. The manifold $W\cup X_0\cup X_1$ is our $X$.

Let $a_i, b_i$ be the generators for the 1st homology (over $\ZZ$) of $T\subset Y_i$. Let $f_i$ be the generator for the fiber class in $Y_i$. Without loss we can assume the Poincare dual of $\omega_4$ is $a_0$; it follows that $a_0=0$ in $X_0$. In $X_1$, clearly $f_1=0$.
The cobordism $W$ deformation retracts onto the space $Y_0\cup Y_1\cup{\rm core}(H)$. Under the attaching of ${\rm core}(H)$ to $Y_0\cup Y_1$, $f_0=f_1$. Therefore the generators of $H_1(X)$ are $b_0, a_1, b_1$. It is  seen that these are also the generators for $H_1(X\cup X)$ of the doubled manifold, so $b_1(X\cup X)=3$.

Let $T_0$ be the torus in $Y_0$ that has intersection number $\pm 1$ with the generator $a_0$. Let $T_1$ be the torus $T\times\mbox{\rm pt}\subset Y_1$. Then $T_0, T_1$ represent generators for $H_2(X)$. There is a disk $D_0$ with $\partial D_0\subset \Sigma_2\times S^1$ such that $D_0\cdot T_0=1$. There is another disk $D_1$ with $D_1\cdot T_1=1$. In the doubled manifold $X\cup X$ the boundaries of these disks match up to form 2-spheres,  call them $S_0, S_1$. If we call $T_0, T_1$ the tori in one copy of $X$, then call $T'_0, T'_1$ the copies in the other. Then the surfaces $T_0, S_0, T'_0, T_1, S_1, T'_1$ are representatives for $H_2(X\cup X)$. The intersection matrix consists of two diagonal blocks. With respect to ordered basis $\{T_0, S_0, T'_0\}$ the first block takes the form
\begin{displaymath}
\left(
\begin{array}{ccc}
  0&1   &0   \\
  1&0   &1   \\
  0&1   & 0  
\end{array}
\right).
\end{displaymath}
The second block is identical. It is verified that each block has maximal positive subspaces of dimension 1, so $b^+_2(X\cup X)=2$.

Return now to equation (\ref{indexE}). A bundle over the doubled manifold that agrees with $E_{\varrho}$ over $X$ can be chosen to have 1st Chern class Poincare dual to $S_0$, because $S_0\cdot T_0=S_0\cdot T'_0=1$, i.e non-zero mod 2. (Other choices such as $S_0+T_0$ etc. are also possible and  will give equivalent results.)

Thus 
$c_1^2=0$ and (\ref{indexE}) evaluates as
\begin{equation}
2\mbox{Ind}E_{\varrho} \equiv 2(4c_2-0)-3(1-3+2)\equiv 0\bmod 4,
\end{equation}
and therefore $\mbox{Ind}E_{\varrho} \equiv  0\bmod 2$. From the preceding considerations it is also straighforward to see that $b_1(X)=3$ and $b^+_2(X)=0$. So
\begin{displaymath}
\nu[\varrho]=\mbox{\rm Ind}E_{\varrho}-3(b_1(X)-b^+_2(X))\equiv 1\bmod 2,
\end{displaymath}
and the lemma follows. \qed

\paragraph{Skein Relations}

Let $J$ denote either an oriented knot or link in the 3-sphere. Let $J_+$, $J_{-}$ and $J_0$ denote in the usual way the knots or links that in a projection differ in a neighborhood of a single crossing. We use the conventions for $J_+$ and $J_{-}$ in \cite[Sect.\ 3]{BD}. If $J_{\pm}$ is modeled on the $x$-$y$-axes  in the plane (the strands agreeing with the standard orientation of the axes), then in $J_+$ the $y$-axis strand goes over the $x$-axis strand.

If $J_{\pm}$ are knots then $J_0$ is a 2-component link. If $J_{\pm}$ are 2-component links and the crossing is between different components of $J_{\pm}$ then $J_0$ is a knot. We have a version of the following well-known theorem but applied only to knots and 2-component links.

\begin{thrm}\label{skein}
Let $J$ denote either  an oriented knot or oriented 2-component link. 
Let  $\Delta_J(t)$ be a finite Laurent polynomial in powers of $t^{1/2}$ that is an oriented isotopy invariant of $J$. Assume 
(1) $\Delta_J(t)=1$ for  the unknot
(2) $\Delta_J(t)=0$ for split 2-component links
and (3) the following Skein rule holds:
\begin{displaymath}
\Delta_{J_+}(t) - \Delta_{J_{-}}(t) = (t^{1/2}-t^{-1/2})\Delta_{J_0}(t),
\end{displaymath}
where if $J$ is a 2 component link then the crossing change is between different strands. 
Then for knots $K$, $\Delta_K(t)$ is exactly the symmetrized and normalized Alexander polynomial of $K$, i.e such that $\Delta_K(t^{-1})=\Delta_K(t)$, and $\Delta_K(1)=1$.
\end{thrm}

\proof see \cite[Proof of Theorem 1.5]{FS}; however note that our $J_{\pm}$ is their $J_{\mp}$. \qed

\section{Proofs}

\paragraph{Theorem~\ref{mainthrm}}
The strategy is to show that for oriented knots  $K$,  $P_K(t)=P_{\omega_0}(K^T,\widehat{\Sigma})(t)$ satisfies the conditions of Theorem~\ref{skein} by applying Floer's exact triangle (\cite{BD} will be our general reference). However this will require, for oriented 2-component links $L$ a slight alteration of $P_{\omega_1}(L^C,\widehat{\Sigma})(t)$ before the Skein relation is satisfied.

We let
\begin{displaymath}
\widehat{P}_J(t) =
\left\{
\begin{array}{ll}
P_{\omega_0}(J^T,\widehat{\Sigma})(t) &\mbox{\rm if $J$  is a knot}\\
-(t^{1/2}-t^{-1/2})P_{\omega_1}(J^C,\widehat{\Sigma})(t) & \mbox{\rm if $J$ is a link}.
\end{array}\right.
\end{displaymath}

According to Floer, for oriented knots $K$ there is an exact sequence of the form
\begin{displaymath}
\longrightarrow
\widetilde{I}_i(K^T_+)_{\omega_0}
\stackrel{a}{\longrightarrow}
\widetilde{I}_i(K^T_{-})_{\omega_0}
\stackrel{b}{\longrightarrow}
\widetilde{I}_i(K^{\Sigma_2}_0)_{\omega_2}
\stackrel{c}{\longrightarrow}
\widetilde{I}_{i-1}(K^T_+)_{\omega_0}
\longrightarrow
\end{displaymath}
where the maps are induced by the various surgery cobordisms, essentially adding 2-handles to the product cobordism along a knot. 
The maps $a$, $b$ are degree zero and $c$ is degree $-1$ applying Lemma~\ref{degree} to the surgery cobordisms $W_a$, $W_b$ and $W_c$ respectively. $W_a$ and $W_b$ are obtained by adding two handles along homologically trivial knots with $-1$-framing. As such $b_1$ is the same as the ends of the cobordism.  Clearly for $W_a$ and $W_b$, $b_2^+=0$. By Lemma~\ref{degree} this gives a dimension shift of $0\bmod 2$. On the other hand, $W_c$ is obtained by adding a 2-handle along a knot which is a generator for 1st homology, so $b_1$ drops by 1 in $W_c$, however $b^+_2=0$. This results in a dimension shift of $-1\bmod 2$.

Since all the surgery cobordisms are connected, the action of $\mu(x_0)$ in each group commutes with $a$, $b$ and $c$. A similar statement is true for the action of $\mu(\widehat{\Sigma})$ on each group, because all the various oriented surfaces $\widehat{\Sigma}$ are homologous in the surgery cobordisms.

Therefore the above exact sequence respects the decomposition of each group into the $\pm2 $ eigenspaces of $\mu(x_0)$ and the generalized eigenspaces for $\mu(\widehat{\Sigma})$ contained therein. Thus
\begin{eqnarray*}\label{s-rel}
P_{\omega_0}(K^T_+,\widehat{\Sigma})
-P_{\omega_0}(K^T_{-},\widehat{\Sigma})
&=&-P_{\omega_2}(K_0^{\Sigma_2},\widehat{\Sigma})\\
&=&-(t^{-1}-2+t)P_{\omega_1}(K_0^C,\widehat{\Sigma}).
\end{eqnarray*}
where the last equality is obtained by Lemmas~\ref{product} and \ref{sigma2}. It immediately follows that 
\begin{equation}
\widehat{P}_{K_+}(t) - \widehat{P}_{K_{-}}(t) 
= (t^{1/2}-t^{-1/2}) \widehat{P}_{K_0}(t).
\end{equation}
In the case of a 2-component oriented link $L$ where the crossing change is between different components of the link, we have the exact sequence
\begin{displaymath}
\longrightarrow
\widetilde{I}_i(L^C_+)_{\omega_1}
\stackrel{a}{\longrightarrow}
\widetilde{I}_i(L^C_{-})_{\omega_1}
\stackrel{b}{\longrightarrow}
\widetilde{I}_i(L^{T}_0)_{\omega_0}
\stackrel{c}{\longrightarrow}
\widetilde{I}_{i-1}(L^C_+)_{\omega_1}
\longrightarrow
\end{displaymath}
A similar argument gives
\begin{displaymath}
P_{\omega_1}(L^C_+,\widehat{\Sigma})
-P_{\omega_1}(L^C_{-},\widehat{\Sigma})
=-P_{\omega_0}(L^T_0,\widehat{\Sigma})
\end{displaymath}
and again relation (\ref{s-rel}) holds. For split 2-component links $L'$ the instanton homology must vanish because $\omega_1$ is non-trivial on the 2-sphere $S$ that separates the components of the link; there are no flat $SO(3)$-connections over $S$. Thus 
$0=P_{\omega_1}(L')(t)=\widehat{P}_{L'}(t)$. For the unknot $U$, $U^T$ is clearly $T\times S^1$ so $\widehat{P}_U(t)=P_{\omega_0}(U^T)(t)=1$ by Lemma~\ref{toruscal}. This completes the proof. \qed

\paragraph{Corollary~\ref{maincor}}

In our definition of $\omega_0$ (defined for $K^T$) we could have also twisted it by the complex line that has Chern class Poincare dual to an oriented fiber $\{\mbox{\rm pt}\}\times S^1$ in the subset 
$({T-D^{\circ}})\times S^1$. We could made corresponding changes for $\omega_1$ and  $\omega_2$ (recall these defined for $L^C$ and $L^{\Sigma_2}$  respectively), ensuring that these are compatible over the knot or link complement with each other. We will not need to change $\omega_3$ (defined for $\Sigma_2\times S^1$) but will need a similar change to $\omega_4$ (defined for $T\times S^1$). Denote the changed complex lines by $\omega'_0$, $\omega_1'$, 
$\omega_2'$ and $\omega_4'$ respectively.

Lemma~\ref{product} remains unchanged with $\omega_1$ and $\omega_2$ replaced with $\omega_1'$ and $\omega_2'$ respectively. Excision does not require $\omega_3$ to be changed.

Lemma~\ref{toruscal} remains the same with $\omega_4$ replaced by $\omega_4'$. There is a diffeomorphism that moves $\omega_4'$ back to $\omega_4$.

Then the proof of Theorem~\ref{mainthrm} goes through with $\omega_i$ replaced with $\omega'_i$. In particular the conclusion about $\widehat{P}_J(t)$, defined using $\omega_i'$, remains the same. Therefore we now assume the $\omega_i$'s are changed to $\omega_i'$'s.

Then, according to Floer (see \cite{BD}) for oriented knots $K$ there is an exact sequence of the form
\begin{displaymath}
\longrightarrow
\widetilde{I}_i(K^D_+)_{\omega'}
\stackrel{a}{\longrightarrow}
\widetilde{I}_i(K^D_{-})_{\omega'}
\stackrel{b}{\longrightarrow}
\widetilde{I}_i(K^{C}_0)_{\omega_1'}
\stackrel{c}{\longrightarrow}
\widetilde{I}_{i-1}(K^D_+)_{\omega'}
\longrightarrow
\end{displaymath}
Repeating the argument of the proof of Theorem~\ref{mainthrm} we have
\begin{displaymath}
Q_{K_+}(t) - Q_{K_{-}}(t) = -P(K_0^C,\widehat{\Sigma})_{\omega_1'}(t),
\end{displaymath}
and therefore
\begin{displaymath}
(t-2+t^{-1})(Q_{K_+}(t) - Q_{K_{-}}(t)) = (t^{1/2}-t^{-1/2})\Delta_{K_0}(t).
\end{displaymath}
Here we make use of the proof of Theorem~\ref{mainthrm} (under the new $\omega_i'$) that identifies $-(t^{1/2}-t^{-1/2})P(K_0^C,\widehat{\Sigma})_{\omega_1'}(t)$ with $\Delta_{K_0}(t)$.
Thus  there is a universal constant $C$ such that
\begin{displaymath}
(t-2+t^{-1})Q_K(t) +C = \Delta_K(t),
\end{displaymath}
where $\Delta_K(t)$ is the symmetrized and normalized Alexander polynomial. By considering the unknot, which has $Q_K(t)=0$, we must have $C=1$. The corollary follows immediately. \qed

\section*{Contact}

Yuhan Lim \newline
Murphys, CA\newline
ylim583@yahoo.com

\end{document}